\documentclass[12pt,a4paper]{article}
\usepackage[ngerman]{babel}
\usepackage[margin=4cm]{geometry}
\usepackage[utf8]{inputenc}
\usepackage{avant}
\usepackage{graphicx}
\usepackage{pslatex}
\usepackage{textcomp}
\usepackage{gensymb}
\usepackage{setspace}
\usepackage{lmodern}
\usepackage{amsmath,amssymb, amsthm}
\usepackage{cases}

\newlength{\bibitemsep}
\newlength{\bibparskip}
\let\oldthebibliography\thebibliography
\renewcommand\thebibliography[1]{%
  \oldthebibliography{#1}%
  \setlength{\parskip}{\bibitemsep}%
  \setlength{\itemsep}{\bibparskip}%
}

\DeclareRobustCommand{\stirling1}{\genfrac{[}{]}{0pt}{}}

\numberwithin{equation}{section}
\newtheorem{theorem}{Theorem} [section]
\newtheorem{lemma}[theorem]{Lemma} 
\newtheorem{ex}[theorem]{Example} 
\newtheorem{cor}[theorem]{Corollary}
\begin{document}
\pagenumbering{arabic}

\title{Moment Generating Stirling Numbers of the first kind and Applications}
\author{Ludwig Frank 
\\ {\small Faculty of Computer Science, Technical University of Applied Sciences,} 
\\ {\small Hochschulstraße 1, D-83024 Rosenheim, Germany,}
\\ {\small ludwig.frank@th-rosenheim.de}
}
\date{}
\maketitle
\begin{quote}
\begin{small}
\textbf{Abstract.} 
In this paper we introduce and investigate moment generating Stirling numbers of the first kind, "`MSN1"'. They are inverses of MSN2's, which make the representation of the moments for a lot of statistical distributions in closed formulas possible. Both MSN1's and MSN2's are related to the r-Stirling numbers, and extend their properties to any real third parameter. If the third parameter is a nonnegative integer, r-Stirling numbers can be converted to MSN's and vice versa. 
\end{small}



\end{quote}

\section{Introduction}\  

In \cite{Fra3} moment generating Stirling numbers of the second kind were defined, "`MSN2"'. They make the representation of the moments for a lot of statistical distributions in closed formulas possible. In the present paper we introduce the moment generating Stirling numbers of the first kind, "`MSN1", and investigate their properties. In some way, they are inverses of MSN2's. Both MSN1's and MSN2's are related to the r-Stirling numbers, (cf. \cite{Bro}) and extend their properties to any real third parameter. If the third parameter is a nonnegative integer, r-Stirling numbers can be converted to MSN's and vice versa. While r-Stirling numbers are useful for combinatorial applications, MSNs allow simpler representations of formulas in statistical applications. 

In Section 2 we present the definition of MSN1's, their basic properties and their interaction with the MSN2's. Finally, Section 3 shows how central, factorial and ordinary moments can be converted into one another using MSN1's (and MSN2's). If the moments of a statistical distribution can be expressed by means of MSN2's, MSN1's make it possible to reconstruct their parameters. 

\section{Moment Generating Stirling Numbers of the first kind}
\subsection{Definitions and Basic Concepts} 

The purpose of this section is to show some simple but basic properties of MSN1's and to present the ordinary and exponential generating function, which allow to derive further properties and identities. They show the relation to ordinary SN1's. 

If $s_{i,j}$ are the (ordinary) Stirling numbers of the first kind (with sign) (cf. e.g. \cite{Abramowitz}), we define for integers $i \ge 0$, $j \ge 0$ and real $k$ the \textbf{\textit{moment generating Stirling numbers (of the first kind) $c_{i,j,k}$, "`MSN1"'}}, by setting
\begin{equation}
	c_{i,j,k}=\sum_{r=j}^{i} \binom{r}{j} \cdot (-k)^{r-j} \cdot s_{i,r}.
	\label{def}
\end{equation}
For $i <0$ or $j <0$ we set $c_{i,j,k}=0$.

The following lemma summarizes some obvious and often used properties:
\begin{lemma} \label{lemma1}
For integers $i \ge 0$, $j \ge 0$ and real $k$, the following properties hold:
\begin{align}
  \text{a) }  &c_{i,j,0}=s_{i,j},  \label{a2}\\
  \text{b) }  &c_{i,j,k}=0, \text { if } j>i, \label{a3}\\
	\text{c) }  &c_{i,i,k}=1, \label{a4}\\
	\text{d) }  &c_{i+1,i,k}=s_{i+1,i}- k \cdot (i+1)=-\binom{i+1}{2}-k \cdot (i+1)
\end{align}
\end{lemma}
\textit{Proof.} a) and b) follow from the definition, c) holds due to $s_{i,i}=1$ (cf. \cite{Abramowitz}) and d) follows from (\ref{def}) and $s_{i+1,i}=-\binom{i+1}{2}$ (cf. e.g. \cite{Abramowitz}). $\qed$

Now, we turn to the ordinary generating function, which is closely related to the generating function of the ordinary Stirling numbers. With 
$(x)_i=x\cdot(x-1)\cdot \cdots \cdot(x-i+1)$, $i \ge 0$, we know from \cite{Abramowitz} that the ordinary generating function of $s_{i,j}$ is 
\begin{equation}
F_i(x)=\sum_{r=0}^{i}s_{i,r}\cdot x^r=(x)_i. \label{def_gen_f_sn}
\end{equation}
For the generating function 
\begin{equation}
G_{i,k}(x)=\sum_{r=0}^{i}c_{i,r,k}\cdot x^r
\end{equation} 
we obtain immediately
\begin {theorem} \ 
For integers $i \ge 0$, $j \ge 0$, $n \ge 0$ and real $k$, the following properties hold:
\begin{align}
&\text{a)  } G_{i,k}(x)=F_{i}(x-k)=(x-k)_i, \label{a8}\\
&\text{b)  } G_{i+n,k}(x)=G_{n,k}(x) \cdot G_{i,k+n}(x),\label{a9} \\
&\text{c)  } G_{i+1,k}(x)=G_{1,k}(x) \cdot G_{i,k+1}(x), \label{a12} \\
&\text{d)  } G_{i,k}(x)\cdot(x-k-i)=G_{i,k+1}(x)\cdot (x-k), \label{a13} \\
&\text{e)  } G_{i,k}(x)\cdot F_k(x)=G_{k,i}(x)\cdot F_i(x), \label{a14} \\
&\text{f)  } F_{i+k}(x)=F_{k}(x)\cdot G_{i,k}(x), \text{ if } k \ge 0.
\end{align}
\end{theorem}
\textit{Proof.} a) 
\begin{align*}
&G_{i,k}(x)=\sum_{j=0}^{i}c_{i,j,k}\cdot x^j= 
\sum_{j=0}^{i}\sum_{r=j}^{i} \binom{r}{j} \cdot (-k)^{r-j} \cdot s_{i,r} \cdot x^j=\\
&\sum_{r=0}^{i}s_{i,r} \cdot\sum_{j=0}^{r}\binom{r}{j} \cdot (-k)^{r-j} \cdot x^j=
\sum_{r=0}^{i}s_{i,r} \cdot (x-k)^r=F_i(x-k)=(x-k)_{i}.
\end{align*}
b) Due to a)
\begin{align*}
&G_{i+n,k}(x)= (x-k)\cdot (x-k-1)\cdot \ldots \cdot (x-k-i-n+1)=\\
&(x-k)\cdot \ldots \cdot (x-k-n+1)\cdot (x-k-n)\cdot \ldots \cdot (x-k-i-n+1)=\\
&G_{n,k}(x)\cdot G_{i,k+n}(x).
\end{align*}
c) follows from b) by setting $n=1$. \\
d) Due to a)
\begin{align*}
&G_{i,k}(x)\cdot (x-k-i)=(x-k)\cdot \ldots \cdot (x-k-i)=\\
&(x-k) \cdot G_{i,k+1}(x).
\end{align*}
e) Due to b)
\begin{align*}
&F_k(x) \cdot G_{i,k}(x)=
x\cdot \ldots \cdot (x-k+1)\cdot(x-k)\cdot \ldots \cdot (x-k-i+1)=\\
&x\cdot \ldots \cdot (x-i+1)\cdot (x-i)\cdot \ldots \cdot (x-k-i+1)=
F_i(x) \cdot G_{k,i}(x).
\end{align*}
f) follows from b) for $k=0$:
\begin{align*}
F_{i+n}(x)=G_{i+n,0}(x)=G_{n,0}(x)\cdot G_{i,n}(x)=F_{n}(x)\cdot G_{i,n}(x).
\end{align*}
Set $n=k$, then $F_{i+k}(x)=F_{k}(x)\cdot G_{i,k}(x)$. $\qed$

If $S_{i,r}$ are the (ordinary) Stirling numbers of the second kind, then  $x^i =\sum_{r=0}^{i}S_{i,r} \cdot (x)_r$ (cf. e.g. \cite{Abramowitz}). The following theorem follows from this fact:
\begin {theorem} \ 
\begin{align}
(x-k)^i=\sum_{r=0}^{i} S_{i,r}\cdot G_{r,k}(x)
\end{align}
\end{theorem}

Now we turn to the exponential generating function
\begin{align}
E_{j,k}(x)=\sum_{i=j}^{\infty} c_{i,j,k} \cdot \frac{x^i}{i!}. \label{exp_gen_1}
\end{align}
\begin {theorem} \ For $|x|<1$
\begin{align}
E_{j,k}(x)=\left( \frac{1}{1+x}\right )^k\cdot \frac{\left(\ln(1+x)\right)^j}{j!}. \label{exp_gen_2}
\end{align}
\end{theorem}
\textit{Proof.} With the exponential generating function of the (ordinary) Stirling numbers of the first kind $s_{i,j}$
\begin{align*}
E_{j,0}(x)=\sum_{i=j}^{\infty} c_{i,j,0} \cdot \frac{x^i}{i!}=
\sum_{i=j}^{\infty} s_{i,j} \cdot \frac{x^i}{i!}=\frac{\left(\ln(1+x)\right)^j}{j!}
\end{align*}
(cf. e.g. \cite{Abramowitz} and Definition (\ref{def})) we obtain
\begin{align*}
&E_{j,k}(x)=\sum_{i=j}^{\infty} c_{i,j,k} \cdot \frac{x^i}{i!}=
\sum_{i=j}^{\infty} \sum_{r=j}^{i} \binom{r}{j}\cdot (-k)^{r-j} \cdot s_{i,r} \cdot \frac{x^i}{i!}=\\
&\sum_{r=j}^{\infty} \binom{r}{j}\cdot (-k)^{r-j}\sum_{i=r}^{\infty} s_{i,r} \cdot \frac{x^i}{i!}
=\sum_{r=j}^{\infty} \binom{r}{j}\cdot (-k)^{r-j}\cdot \frac{1}{r!}\cdot \left( \ln(1+x)\right)^r=\\
&\frac{1}{j!}\cdot \sum_{r=j}^{\infty} (-k)^{r-j} \cdot \frac{\left( \ln(1+x)\right)^{r}}{(r-j)!}=
\frac{1}{j!}\cdot \sum_{s=0}^{\infty} (-k)^{s} \cdot \frac{\left( \ln(1+x)\right)^{s+j}}{s!}=\\
&\frac{1}{j!}\cdot \left( \ln(1+x)\right)^j \cdot \exp \left( -k \cdot \ln(1+x)\right)=
\frac{1}{j!}\cdot \left( \ln(1+x)\right)^j \cdot \frac{1}{(1+x)^k}.\qed
\end{align*} 

\subsection{Properties and Identities of MSN1's}
Because the computation of the values $c_{i,j,k}$ according the definition is boring and error prone, the current section presents some formulae to compute them. The main aids for that are the formulae 
\begin{align}
&c_{i,j,k}=\frac{G_{i,k}^{(j)}(0)}{j!}, \label{c_ijk_gen_f_1}\\ 
&(f\cdot g)^{(j)}=\sum_{r=0}^{j} \binom{j}{r} \cdot f^{(r)} \cdot g^{(j-r)} \label{c_ijk_gen_f_2}\\ 
&\text{and}  \nonumber \\ 
&(f_1 \cdot \ldots \cdot f_n)'=
(f_1 \cdot \ldots \cdot f_n)\cdot \left(\frac{f_1'}{f_1} + \ldots + \frac{f_n'}{f_n} \right). \label{c_ijk_gen_f_3}
\end{align}
for differentiable functions $f$, $f_1, \ldots, f_n$ and $g$.
\begin {theorem} For integers $i \ge 0$, $j \ge 0$ and real $k$
\begin{align}
&\text{a)  } c_{i,0,k}=(-k)_i,\label{c_ijk_1} \\ 
&\text{b)  } c_{i,0,1}=(-1)^i \cdot i!,\label{c_ijk_2} \\ 
&\text{c)  } c_{i,1,k}=(-1)^{i+1} \cdot \binom{k+i-1}{i}\cdot i! \cdot (H_{k+i-1}-H_{k-1}) \label{c_ijk_3}\\ 
&\text{ for nonnegative integer k}, \text { if } H_n=\sum_{r=1}^n \frac{1}{r} \text { are the harmonic numbers,} \nonumber \\ 
&\text{d)  } c_{i,1,1}=(-1)^{i+1} \cdot i! \cdot H_i, \label{c_ijk_4}\\ 
&\text{e)  } c_{i,j-1,1}=s_{i+1,j},  \label{c_ijk_5}\\ 
&\text{f)  } c_{i+1,j,-1}= s_{i,j-1}+s_{i,j}. \label{c_ijk_6}
\end{align}
\end{theorem}
\textit{Proof.} a) Due to (\ref{c_ijk_gen_f_1})
\begin{align*}
c_{i,0,k}=G_{i,k}(0)=F_{i}(-k)=(-k)_i.
\end{align*}
b) follows from a) by $k=1$. \\ 
c) Due to (\ref{c_ijk_gen_f_1}) and (\ref{c_ijk_gen_f_3})
\begin{align*}
c_{i,1,k}=&G'_{i,k}(0)=F'_i(-k)=F_i(-k) \cdot \left(-\frac{1}{k}- \ldots -\frac{1}{k+i-1} \right)=\\
&(-1)^{i+1} \cdot \frac{(k+i-1)!}{(k-1)!} \cdot \left(\frac{1}{k}+ \ldots +\frac{1}{k+i-1} \right)=\\
&\binom{k+i-1}{i}\cdot i! \cdot (H_{k+i-1}-H_{k-1})
\end{align*}
d) follows from c) for $k=1$. \\ 
e) Due to (\ref{a13}) with $k=0$ $G_{i,1}(x)\cdot x=G_{i,0}(x)\cdot(x-i)$. Differentiating this equation $j$ times results in  
\begin{align*}
\sum_{r=0}^1 \binom{j}{r} \cdot G_{i,1}^{(j-r)}(0) \cdot \left. x^{(r)} \right|_{x=0} =
\sum_{r=0}^1 \binom{j}{r} \cdot G_{i,0}^{(j-r)}(0) \left.\cdot (x-i)^{(r)}\right|_{x=0} ,
\end{align*}
from which we obtain
\begin{align*}
&G_{i,1}^{(j)}(0) \left.\cdot x ^{(0)}\right|_{x=0} +j \cdot \left. G_{i,1}^{(j-1)}(0)   \cdot (x)^{(1)} \right|_{x=0} =\\
&G_{i,0}^{(j)}(0) \cdot (-i) + j \cdot G_{i,0}^{(j-1)}(0)\cdot 1,
\end{align*}
which implies
\begin{align*}
j \cdot (j-1)! \cdot c_{i,j-1,1} = 
s_{i,j}\cdot j! \cdot (-i)+ j \cdot (j-1)! \cdot s_{i,j-1}
\end{align*}
and so
\begin{align*}
c_{i,j-1,1}=-i \cdot s_{i,j}+s_{i,j-1}=s_{i+1,j}. 
\end{align*} \\ 
f) Due to (\ref{a12}) with $k=-1$
\begin{align*}
G_{i+1, -1}(x)=G_{1,-1}(x) \cdot G_{i,0}(x) =(x+1)\cdot G_{i,0}(x)
\end{align*}
and due to (\ref{c_ijk_gen_f_1}) and (\ref{c_ijk_gen_f_2}) by differentiating this equation $j$ times
\begin{align*}
&c_{i+1,j,-1}=
\frac{1}{j!} \cdot \left(\sum_{r=0}^j \binom{j}{r} \cdot G_{1,-1}^{(r)}(0)\cdot G_{i,0}^{(j-r)}(0)\right) =\\
&\frac{1}{j!} \cdot \left(G_{1,-1}^{(0)}(0)\cdot G_{i,0}^{(j)}(0)+ j \cdot
 G_{1,-1}^{(1)}(0)\cdot G_{i,0}^{(j-1)}(0)\right)= \\
&\frac{1}{j!} \cdot \left(1 \cdot j! \cdot c_{i,j,0} + 
j \cdot 1 \cdot (j-1)! \cdot c_{i,j-1,0}\right) =
s_{i,j} +s_{i,j-1}. \qed
\end{align*}

\begin {theorem} For integer $i \ge 0$ and real $k_1$ and $k_2$
\begin{align}
\sum_{r=0}^{i} c_{i,r,k_1} \cdot k_2^r=(k_1-k_2)_i. \label{sum_k1_k2}
\end{align}
\end{theorem}
\textit{Proof.} Due to (\ref{a8})
\begin{align*} 
\sum_{r=0}^{i} c_{i,r,k_1} \cdot k_2^r=G_{i,k_1}(k_2)=F_i(k_2-k_1)=(k_2-k_1)_i. \qed
\end{align*}

\begin{cor} For integer $i \ge 0$ and real $j$, $k$, $k_1$ and $k_2$
\begin{align}
&\text{ a) } \sum_{r=0}^{i} c_{i,r,k} = c_{i,0,k-1},\\
&\text{ b) } \sum_{r=0}^{i} c_{i,r,2} = (-1)^i \cdot i!,\\
&\text{ c) } \sum_{r=0}^{i} c_{i,r,k} \cdot (j+k)^r= (j)_i,\\
&\text{ d) } \sum_{r=0}^{i} c_{i,r,k} \cdot (i+k)^r=i!,\\
&\text{ e) } \sum_{r=0}^{i} c_{i,r,k} \cdot m^r =0 
\text{ iff } m \in \{k, \ldots, k+i-1\} \text{ for } i>0,\\
&\text{ f) } \sum_{r=0}^{i} c_{i,r,k} =0
\text{ iff } k \in \{1,0,-1, \ldots, -i+2\} \text{ for }i>0,\\
&\text{ g) } c_{i,0,2}=(-1)^i \cdot (i+1)!. 
\end{align}
\end {cor}
\textit{Proof.} a) Due to (\ref{sum_k1_k2}) and (\ref{c_ijk_1}) 
$\sum_{r=0}^{i} c_{i,r,k}=(1-k)_i=c_{i,0,k-1}$. b) follows from a) due to (\ref{c_ijk_2}). c) The left side is $G_{i,k}(j+k)=F_i(j+k-k)=(j)_i$ due do (\ref{def_gen_f_sn}). d) The left side is $G_{i,k}(i+k)=F_i(i)=i!$. e) The left side is $G_{i,k}(m)=(m-k)_i=(m-k)\cdot(m-k-1)\cdot \ldots \cdot(m-k-i+1)$. f) The left side is $G_{i,k}(1)=(1-k)_i=(1-k)\cdot(1-k-1)\cdot \ldots \cdot(1-k-i+1)$. g) The left side is $G_{i,2}(0)=F_i(-2)=(-2)\cdot(-3)\cdot \ldots \cdot(-2-i+1)=
(-1)^i \cdot (i+1)! \qed$.

More complicated computations of $c_{i,j,k}$ can be made easier by the following recursions. They are similar to those of the MSN2's in \cite{Fra3} and the ordinary Stirling numbers of the first and second kind (cf. \cite{Abramowitz} and \cite{DLMF}).

\begin {theorem} For real $k$ and integers $n \ge 0$, $i \ge 0$ and $j \ge 0$
\begin{align}
&\text{ a) } c_{i+n,j,k} = \sum_{r=0}^{j} c_{n,r,k}\cdot c_{i,j-r,k+n}, \\
&\text{ b) } c_{i+1,j,k} = c_{i,j-1,k} - (i+k)\cdot c_{i,j,k},\label{rel_to_r}\\ 
&\text{ c) } c_{i+1,j,k} = c_{i,j-1,k+1} - k \cdot c_{i,j,k+1}.
\end{align}
\end{theorem} 
\textit{Proof.} a) From $G_{i+n,k}(x)=G_{n,k}(x) \cdot G_{i,k+n}(x)$, (cf. (\ref{a9})), it follows by differentiation $j$-times and due to (\ref{c_ijk_gen_f_2}) 
\begin{align*}
c_{i+n,j,k} =&\frac{1}{j!} \cdot G_{i+n}^{(j)}(0)=\frac{1}{j!} \cdot 
\sum_{r=0}^j \binom{j}{r} \cdot G_{n,k}^{(r)}(0) \cdot G_{i,k+n}^{(j-r)}(0)= \\
&\frac{1}{j!} \cdot \sum_{r=0}^j \binom{j}{r} \cdot r! \cdot c_{n,r,k} \cdot
(j-r)! \cdot c_{i,j-r,k+n}=\\
& \sum_{r=0}^j c_{n,r,k} \cdot c_{i,j-r,k+n}.
\end{align*}
b) From $G_{i+1,k}(x)=G_{i,k}(x)\cdot G_{1,k}(x)=G_{i,k}(x)\cdot (x-i-k)$  follows by differentiating $j$-times  
\begin{align*}
c_{i+1,j,k} =&\frac{1}{j!} \cdot G_{i+1}^{(j)}(0)=
\frac{1}{j!} \cdot \sum_{r=0}^j \binom{j}{r} \left. \cdot (x-i-k)^{(r)}\right|_{x=0}\cdot G_{i,k}^{(j-r)}(0)=\\
&\frac{1}{j!} \cdot \left(-(i+k)\cdot j! \cdot c_{i,j,k}+ j \cdot (j-1)!\cdot c_{i,j-1,k} \right)=\\
&-(i+k)\cdot c_{i,j,k}+c_{i,j-1,k}.
\end{align*}
c) By differentiating 
$G_{i+1,k}(x)=G_{1,k}(x) \cdot G_{i,k+1}(x)= (x-k)\cdot G_{i,k+1}(x)$ (cf. (\ref{a12})) we obtain
\begin{align*}
&c_{i+1,j,k} =\frac{1}{j!} \cdot \sum_{r=0}^j \binom{j}{r} \left. \cdot (x-k)^{(r)}\right|_{x=0}\cdot G_{i,k+1}^{(j-r)}(0)=\\
&\frac{1}{j!} \cdot \left(-k \cdot j! \cdot c_{i,j,k+1}+j \cdot 
(j-1)!\cdot c_{i,j-1,k+1}\right)=\\
&-k \cdot c_{i,j,k+1}+c_{i,j-1,k+1}. \qed
\end{align*}

By means of (\ref{rel_to_r}), the relationship to the r-Stirling numbers of the first kind $\stirling1{i}{j}_r$ can be shown.
\begin {theorem} For integers $i \ge 0$, $j \ge 0$ and $r \ge 0$
\begin{align}
&\stirling1{i}{j}_r=(-1)^{i-j}\cdot c_{i-r,j-r,r}, \label{r-stir-1}\\
&\stirling1{i}{j}_r=0 \text{ for } i<r, \label{r-stir-2}\\
&\stirling1{r}{j}_r=\delta_{j,r}. \label{r-stir-3}
\end{align}
\end{theorem} 
\textit{Proof.} It is shown that with the settings (\ref{r-stir-1}), (\ref{r-stir-2}) and (\ref{r-stir-3}) the defining recursion and the initial conditions of the r-Stirling numbers of the first kind are fulfilled. 
Due to (\ref{rel_to_r})
\begin{align*}
\stirling1{i}{j}_r=&(-1)^{i-j}\cdot c_{i-r,j-r,r} = \\
&(-1)^{i-j} \cdot \left(c_{i-r-1,j-r-1,r}-(i-r+r-1)\cdot c_{i-r-1,j-r,r}\right)=\\
&(-1)^{i-j} \cdot \left(c_{i-r-1,j-r-1,r}-(i-1)\cdot c_{i-r-1,j-r,r}\right)=\\
&(-1)^{i-j} \cdot \left( (-1)^{i-j} \stirling1{i-1}{j-1}_r -
(i-1) \cdot (-1)^{i-j-1} \cdot  \stirling1{i-1}{j}_r\right)=\\
&\stirling1{i-1}{j-1}_r+(i-1) \cdot \stirling1{i-1}{j}_r.
\end{align*}
This is the defining recursion in \cite{Bro}. If $i<r$, then 
\begin{align*}
\stirling1{i}{j}_r=(-1)^{i-j} \cdot c_{i-r,j-r,r}=0
\end{align*}
due to the definition of $c_{i,j,k}$. Further 
\begin{align*}
\stirling1{r}{j}_r=(-1)^{r-j} \cdot c_{r-r,j-r,r}=
(-1)^{r-j} \cdot c_{0,j-r,r}=\delta_{j,r},
\end{align*}
also because of the definition. So, also the initial conditions in \cite{Bro} are fulfilled. Therefore the values $(-1)^{i-j}\cdot c_{i-r,j-r,r}$ are the r-Stirling numbers. $\qed$

The preceding theorem allows to translate all results from the "`r-Stirling number language"' to the "`MSN1 language"' and vice versa. Also the combinatorial interpretation of the r-Stirling numbers can  be transferred to that of the MSN1's.

Now, let's switch to the inverse of the matrix $(c_{i,j,k})_{i,j}$. In \cite{Fra3} it was shown, that $(c_{i,j,k})_{i,j}$ is the right inverse of $\left( \frac{b_{i,j,k}}{j!}\right)_{i,j}$, where $b_{i,j,k}$ are the \textbf{\textit{moment generatig Stirling numbers of the second type, "`MSN2"'}}, defined in \cite{Fra3} by
\begin{equation}
	b_{i,j,k}=\sum_{r=0}^{j} \binom{j}{r} \cdot (-1)^{j-r} \cdot (r+k)^i.
	\label{a1}
\end{equation}
Here, it is proved that it is also the left inverse and in a more general sense the product $(c_{i,j,k_1})_{i,j} \cdot (b_{i,j,k_2})_{i,j}$ is computed.
\begin {theorem} For integers $i \ge 0$, $j \ge 0$ and real $k$  
\begin{align}
&\text{ a) } \sum_{r=j}^{i}c_{i,r,k}\cdot b_{r,j,k}=j! \cdot \delta_{i,j}, \label{inv1} \\
&\text{ b) } \sum_{r=j}^{i}c_{i,r,k_1}\cdot b_{r,j,k_2} =i! \cdot \binom{k_2-k_1}{i-j} \text{ for integers } k_2 \ge k_1.   \label{inv2} 
\end{align}
\end{theorem} 
\textit{Proof.} a) Due to Corollary 2.7 in \cite{Fra3} and the inversion formula for the ordinary Stirling numbers $s_{i,j}$ and $S_{i,j}$ of the first resp. second type
\begin{align*}
&\sum_{r=j}^{i}c_{i,r,k}\cdot b_{r,j,k}=
\sum_{r=j}^{i}\sum_{s=r}^{i} \binom{s}{r} \cdot (-k)^{s-r} \cdot s_{i,s} \cdot b_{r,j,k}=\\
&\sum_{s=j}^{i} s_{i,s} \cdot\sum_{r=j}^{s} \binom{s}{r} \cdot (-k)^{s-r} \cdot b_{r,j,k}=
\sum_{s=j}^{i} s_{i,s} \cdot b_{s,j,0}=\\ 
&\sum_{s=j}^{i} s_{i,s} \cdot S_{s,j} \cdot j!=
j! \cdot \delta_{i,j}.
\end{align*}
b) Due to \cite{Fra3}, Theorem 2.13,
\begin{align*}
&\sum_{r=j}^{i}c_{i,r,k_1}\cdot b_{r,j,k_2} = 
\sum_{r=j}^{i}c_{i,r,k_1}\cdot \sum_{s=0}^{k_2-k_1} \binom{k_2-k_1}{s} \cdot
b_{r,j+s,k_1}=\\
&\sum_{s=0}^{k_2-k_1} \binom{k_2-k_1}{s} \cdot 
\sum_{r=j}^{i}c_{i,r,k_1}\cdot b_{r,j+s,k_1}=
\sum_{s=0}^{k_2-k_1} \binom{k_2-k_1}{s} \cdot \delta_{i,j+s} \cdot (j+s)!=\\
&\binom{k_2-k_1}{i-j} \cdot i!. \qed
\end{align*}

Also a convolution theorem, like the ones for the ordinary Stirling numbers and the MSN2's exists here. 

\begin {theorem} For integers $i \ge 0 $, $j_1 \ge 0 $ and $j_2\ge 0 $ and real $k_1$ and $k_2$
\begin{align}
\sum_{r=0}^{i}\binom{i}{r} \cdot c_{i-r,j_1,k_1}\cdot c_{r,j_2,k_2}=
\binom{j_1+j_2}{j_1} \cdot c_{i,j_1+j_2,k_1+k_2}.
\end{align}
\end{theorem} 
\textit{Proof.} Due to (\ref{exp_gen_2}), for the exponential generating function 
\begin{align*}
E_{j_1,k_1}(x) \cdot E_{j_2,k_2}(x) =
\binom{j_1+j_2}{j_1}\cdot E_{j_1+j_2,k_1+k_2}(x). 
\end{align*}
By differentiating $i$-times, due to $c_{i,j,k}=E_{j,k}^{(i)}(0)$, we obtain the equation
\begin{align*}
\sum_{r=0}^{i}\binom{i}{r}\cdot E_{j_1,k_1}^{(i-r)}(0)\cdot E_{j_2,k_2}^{(r)}(0)=  \binom{j_1+j_2}{j_1} \cdot E_{j_1+j_2,k_1+k_2}^{(i)}(0),
\end{align*}
which implies the assertion. $\qed$

\begin{cor} For integers $i  \ge 0$, $j \ge 0$ and real $k$
\begin{align}
&\text{ a) } \sum_{r=j}^{i} \binom{i}{r} (k)_{i-r} \cdot c_{r,j,k} = s_{i,j},\\
&\text{ b) } \sum_{r=j}^{i} (i)_{r}\cdot (-1)^{i-r} \cdot c_{r,j,k} =c_{i,j,k+1}.
\end{align}
\end {cor}
\textit{Proof.} a) follows from the preceding theorem for $j_1=0$, because $c_{i-r,0,k_1}=(-k_1)_{i-r}$ (cf. (\ref{c_ijk_1})) and so
\begin{align*}
\sum_{r=j}^{i}\binom{i}{r} \cdot (-k_1)_{i-r}\cdot c_{r,j_2,k_2}=
c_{i,j_2,k_1+k_2}.
\end{align*}
Setting $k_1=-k_2$ implies a). Also b) follows from the preceding theorem for $j_1=0$ and $k_1=1$:
\begin{align*}
&c_{i,j,k+1}=\sum_{r=j}^{i} \binom{i}{r} c_{i-r,0,1} \cdot c_{r,j,k}=
\sum_{r=j}^{i} \binom{i}{r} (-1)_{i-r} \cdot c_{r,j,k}=\\
&\sum_{r=j}^{i} \binom{i}{r} \cdot (i-r)!  \cdot(-1)^{i-r} \cdot c_{r,j,k}=
\sum_{r=j}^{i} \frac{i!}{r!} \cdot(-1)^{i-r}\cdot c_{r,j,k}=\\
&\sum_{r=j}^{i} (i)_{r}\cdot (-1)^{i-r} \cdot c_{r,j,k}. \qed
\end{align*}
\begin {theorem} For integers $i  \ge 0$, $j \ge 0$ and real $k_1$ and $k_2$
\begin{align}
\sum_{r=j}^{i} \binom{r}{j} \cdot k_2^{r-j} \cdot c_{i,r,k_1} = c_{i,j,k_1-k_2}
\end{align}
\end{theorem} 
\textit{Proof.} By definition of $c_{i,r,k_1}$
\begin{align*}
&\sum_{r=j}^{i} \binom{r}{j} \cdot k_2^{r-j} \cdot c_{i,r,k_1} = 
\sum_{r=j}^{i} \binom{r}{j} \cdot k_2^{r-j} \cdot 
\sum_{s=r}^{i} \binom{s}{r} \cdot (-k_1)^{s-r} \cdot s_{i,s}= \\
&\sum_{s=j}^{i} s_{i,s} \cdot 
\sum_{r=j}^{s} \binom{s}{r} \cdot \binom{r}{j} \cdot k_2^{r-j} \cdot (-k_1)^{s-r} =\\
&\sum_{s=j}^{i} s_{i,s} \cdot \binom{s}{j} \cdot 
\sum_{r=j}^{s} \binom{s-j}{r-j} \cdot k_2^{r-j} \cdot (-k_1)^{s-r} =\\
&\sum_{s=j}^{i} s_{i,s} \cdot \binom{s}{j} \cdot 
\sum_{t=0}^{s-j} \binom{s-j}{t} \cdot  k_2^{t} \cdot (-k_1)^{s-j-t} =\\
&\sum_{s=j}^{i} s_{i,s} \cdot \binom{s}{j} \cdot (k_2-k_1)^{s-j}=
c_{i,j,k_1-k_2}. \qed
\end{align*}
\begin {cor}
\begin{align}
&\text{ a) } \sum_{r=j}^{i} \binom{r}{j} \cdot k^{r-j} \cdot c_{i,r,k} = 
s_{i,j},\\
&\text{ b) } \sum_{r=j}^{i} \binom{r}{j} \cdot c_{i,r,k} =c_{i,j,k-1},\\
&\text{ c) } \sum_{r=j}^{i} \binom{r}{j} \cdot (-1)^{r-j} \cdot c_{i,r,k} =
c_{i,j,k+1}.
\end{align}
\end {cor}
\textit{Proof.} a) follows from the preceding theorem for $k=k_1=k_1$, because
\begin{align*}
\sum_{r=j}^{i} \binom{r}{j} \cdot k^{r-j} \cdot c_{i,r,k} = c_{i,j,0}= s_{i,j}.
\end{align*}
b) and c) follow by setting $k_2=1$ resp. $k_2=-1$. $\qed$

\section{Applications}\  

In \cite{Fra3} binomial coefficients and the MSN2's allowed to express the central moments
\begin{align}
C_m(X)=M_1\left((X-M_1(X))^m\right)=\sum_{j=0}^m \binom{m}{j} \cdot \left( -M_1(X)\right)^{m-j}  \cdot M_j(X)
\end{align}
of any random variable $X$ on \{0, 1, 2,\ldots\} by the ordinary moments $M_m(X)$ and to express the ordinary moments by the factorial moments
\begin{align}
F_m(X)=M_1(X \cdot (X-1)\cdot \ldots \cdot (X-m+1)) 
\end{align}
as 
\begin{align} 
M_m(X)=
\sum_{j=0}^m b_{m,j,0} \cdot \frac{1}{j!} \cdot F_j(X). \label{FmMm}
\end{align}
Additionally it was shown that
\begin{align} 
C_m(X)=\sum_{j=0}^m b_{m,j,-M_1(X)}\cdot \frac{1}{j!} \cdot F_j(X). \label{CmFm}
\end{align}
The following theorem, which inverts these equations, shows that the MSN1's make it possible to convert any kind of moment to any other one.
\begin {theorem} If $X$ is a random variable on \{0, 1, 2, \ldots\} with existing moments of any order, then 
\begin{align}
&\text{ a) } F_m(X)=\sum_{j=0}^m c_{m,j,0} \cdot  M_j(X),\\
&\text{ b) } F_m(X)=\sum_{j=0}^m c_{m,j,-M_1(X)} \cdot C_j(X),\\
&\text{ c) } M_m(X)=\sum_{j=0}^m \binom{m}{j} \cdot \frac{1}{j!} \cdot \left( M_1(X)\right)^{m-j}  \cdot C_j(X).
\end{align}
\end{theorem} 
\textit{Proof.} a) Due to (\ref{inv1}) and (\ref{FmMm})
\begin{align*}
&\sum_{j=0}^m c_{m,j,0} \cdot  M_j(X)=\sum_{j=0}^m c_{m,j,0} \cdot \sum_{r=0}^j \frac{b_{j,r,0}}{r!} \cdot F_r(X) = \\
&\sum_{r=0}^m \sum_{j=r}^m c_{m,j,0} \cdot \frac{b_{j,r,0}}{r!} \cdot F_r(X)= F_m(X).
\end{align*}
b) Due to (\ref{inv1}) and (\ref{CmFm})
\begin{align*}
&\sum_{j=0}^m c_{m,j,-M_1(X)} \cdot C_j(X)= 
\sum_{j=0}^m c_{m,j,-M_1(X)} \cdot \sum_{r=0}^j b_{j,r,-M_1(X)}\cdot \frac{1}{r!} \cdot F_r(X)=\\
&\sum_{r=0}^m F_r(X) \cdot \frac{1}{r!} \cdot \sum_{j=r}^m c_{m,j,-M_1(X)} \cdot b_{j,r,-M_1(X)}= F_m(X).
\end{align*} 
c) Due to (\ref{FmMm}), b) and Theorem 2.18 in \cite{Fra3}
\begin{align*}
&M_m(X)= \sum_{j=0}^m  \frac{b_{m,j,0}}{j!} \cdot F_j(X)= \sum_{j=0}^m  \frac{b_{m,j,0}}{j!} \cdot \sum_{r=0}^j c_{j,r,-M_1(X)} \cdot C_r(X) = \\
&\sum_{r=0}^m C_r(X) \cdot \sum_{j=r}^m \frac{b_{m,j,0}}{j!}\cdot c_{j,r,-M_1(X)} =
\sum_{r=0}^m C_r(X) \cdot \binom{m}{r} \cdot M_1(X)^{m-r}. \qed
\end{align*}
Just as in \cite{Fra3} the central moments for a lot of probability distributions were possible in closed formulas, here the factorial moments can be determined in simple calculations by means of Theorem 3.1. 
\begin {ex}   
If the random variable $X$ is Poisson distributed with parameter $\lambda$, then 
\begin{align*}
M_m(x)=\sum_{j=0}^m S_{m,j} \cdot \lambda^j=
\sum_{j=0}^m b_{m,j,0} \cdot \frac{\lambda^j}{j!}. 
\end{align*}
By Item a) of the preceding theorem we obtain due to (\ref{inv1})
\begin{align*}
&F_m(X)=\sum_{j=0}^m c_{m,j,0} \cdot M_j(X)=
\sum_{j=0}^m c_{m,j,0} \cdot \sum_{r=0}^j b_{j,r,0} \cdot \frac{\lambda^r}{r!}=\\
&\sum_{r=0}^m \frac{\lambda^r}{r!} \cdot \sum_{j=r}^m c_{m,j,0} \cdot b_{j,r,0}= 
\lambda^m.
\end{align*}
\end{ex}
Of special interest is the phase type distribution $PH(\textbf{\textit{p}},\textbf{\textit{P}})$, with a substochastic matrix $\textbf{\textit{P}}$ and a suitable line vector $\textbf{\textit{p}}$ (cf. \cite{Neuts}, \cite{Neuts_geom} and \cite{Fra2}). In \cite{Fra3}, it was shown, if we assume for $X \sim PH(\textbf{\textit{p}},\textbf{\textit{P}})$ $P(X=0)=0$, (cf. \cite{Fra3}) 
\begin{align}
M_m(X)=\sum_{r=0}^m b_{m,r,1} \cdot \textbf{\textit{p}} \cdot 
\left( \textbf{\textit{P}}\cdot (\textbf{\textit{I-P}}) \right)^{r} \cdot \textbf{\textit{e}},
\end{align}
with a suitable column vector $\textbf{\textit{e}}^T=(1, 1,\ldots ,1)$. This has due to (\ref{inv1}) the consequence 
\begin{align}
F_m(X)=m! \cdot \textbf{\textit{p}} \cdot \textbf{\textit{P}}^{m-1} \cdot (\textbf{\textit{I-P}})^{-m}\cdot \textbf{\textit{e}}, m \ge 1.
\end{align}
This is a system of nonlinear equations for $\textbf{\textit{p}}$ and $\textbf{\textit{P}}$, which allows to reconstruct these values from the moments in principle. 
In a similar way the moments and factorial moments of the recurrence time of a discrete time Markov chain can be expressed in a straight forward way.  

\paragraph{}
Some questions must remain open here, which should be answered in further investigations:
\begin{itemize}
\item What is a combinatorial interpretation of \textit{MSNs} and how can they be applied, if $k$ is not an integer greater or equal zero?
\item How can the parameters of a phase type distribution or a recurrence time reconstructed 
efficiently from the moments.
\end{itemize}

\renewcommand{\refname}{References}
\setlength{\bibitemsep}{0.0\baselineskip}


\begin{thebibliography} {0}
\bibitem{Abramowitz} \textsc{M. Abramowitz and I. A. Stegun} (1970) \textit{Handbook of Mathematical Functions.} National Bureau of Standards, Washington
\bibitem{DLMF} \textsc{D. M. Bressoud} (2015) \textit{NIST Digital Library of Mathematical Functions, Chapter 26, Combinatorial Analysis}. Available online at https://dlmf.nist.gov/26. [Online; accessed 27-May-2021]
\bibitem{Bro}\textsc{A. Z. Broder} (1984) \textit{The r-Stirling numbers.} Discrete Mathematics Vol. 49, 241-259
\bibitem{Fra} \textsc{L. Frank} (2017) \textit{The Binomial and Negative Binomial Distribution in Discrete Time Markov Chains.} Markov Processes and Related Fields 23, 377-400 
\bibitem{Fra2} \textsc{L. Frank} (2022) \textit{The Moments of the Generalized Negative Binomial Distribution and Recurrence Time in Discrete Time Marcov Chains.} Markov Processes and Related Fields 28, 53-85 
\bibitem{Fra3} \textsc{L. Frank} (2022)\textit{Moment Generating Stirling Numbers and Statistical Applications.} arXiv:2206.08837 
\bibitem{Neuts} \textsc{M. F. Neuts} (1978) \textit{Renewal Processes of Phase Type.} Naval research logistics 25, 445-454.
\bibitem{Neuts_geom} \textsc{M. F. Neuts} (1981) \textit{Matrix-Geometric Solutions in Stochastic Models.} Dover Publications Inc. 
\end{thebibliography}
\end{document}